\documentclass[12pt]{amsart}

\usepackage{amsmath}
\usepackage{amsthm}
\usepackage{amsfonts}
\usepackage{amssymb,amsmath,latexsym}
\usepackage{graphicx}

\theoremstyle{plain}
\newtheorem{theorem}{Theorem}[section]
\newtheorem{lemma}{Lemma}[section]

\theoremstyle{definition}

\title[On radial SLE in multiply connected domains]%
{On Radial Stochastic Loewner evolution in multiply connected domains}

\author{Robert~O. Bauer}
\address{Department of Mathematics\\ 
	University of Illinois at Urbana-Champaign\\ 
	1409 West Green Street \\ 
	Urbana, IL 61801, USA}
\email{rbauer@math.uiuc.edu}
\thanks{The research of the first author was supported by 	NSA grant H98230-04-1-0039}
\author{Roland~M. Friedrich}
\address{Institute for Advanced Study\\ 
	Princeton, NJ 08540, USA}
\email{rolandf@ias.edu}
\thanks{The research of the second author was supported by 	NSF grant DMS-0111298.}

\begin{document}

\begin{abstract}
We discuss the extension of radial SLE to multiply connected planar domains. First, we extend Loewner's theory of slit mappings to multiply connected domains by establishing the radial Komatu-Loewner equation, and show that a simple curve from the boundary to the bulk is encoded by a motion on moduli space and a motion on the boundary of the domain. Then, we show that the vector-field describing the motion of the moduli is Lipschitz.  We explain why this implies  that ``consistent," conformally invariant random simple curves are described by multidimensional diffusions, where one component is a motion on the boundary, and the other component is a motion on moduli space. We argue what the exact form of this diffusion is (up to a single real parameter $\kappa$) in order to model boundaries of percolation clusters. Finally, we show that this moduli diffusion leads to random non-self-crossing curves satisfying the locality property if and only if $\kappa=6$. 
\end{abstract}

\maketitle

\section{Introduction}

In this paper we discuss conformally invariant growing random compact sets in multiply connected domains. Our results are meant to provide some of the steps to extend the stochastic Loewner evolution of Schramm, \cite{schramm:2000}, from simply connected domains to multiply connected domains and Riemann surfaces. Based on Loewner's theory of slit mappings, Schramm showed that conformally invariant measures describing random ``simple" curves in simply connected domains, can be encoded into a diffusion on the boundary of the domain if the random simple curves satisfy a ``consistency condition." Furthermore, under an additional but very natural symmetry condition, he showed that the diffusion on the boundary is a multiple of Brownian motion if the random simple curve lives in certain standard domains.

We show that consistent and conformally invariant random simple curves in multiply connected domains can be encoded into a multidimensional diffusion. One component of this diffusion corresponds to a motion on the boundary of the domain and the other components are the moduli of the domain with the random simple curve, grown up to time $t$, removed. The random simple curves we consider in this paper connect the boundary to the bulk (interior) and so the appropriate moduli space is the moduli space of $n$-connected domains with one marked interior point and one marked boundary point. We show that consistency (a Markovian-type property) and conformal invariance essentially determine the diffusion up to the drift of the motion on the boundary. We call such diffusions on moduli space Schiffer diffusions. Under a symmetry condition, familiar from percolation, this drift component can also be identified, leaving a single real parameter $\kappa$. Beginning with the Schiffer diffusion with this drift we show that the resulting family of random growing compact sets satisfies the locality property if $\kappa=6$.

The fundamental observation that diffusion processes on the moduli space of bordered Riemann surfaces with marked points, given its path-wise solutions agree with the geometric constraints, yield the most general way to define probability measures on (simple) curves on surfaces and therefore contain ``ordinary" SLE as a special case, was introduced in~\cite{FK, KBonn}. However, the current article is the first constructive implementation of it, for the radial case and multiply connected domains. For general Riemann surfaces and the chordal case one proceeds along similar lines as described here, but with some necessary modifications. An important role then is played by the so-called ``Hilbert Uniformisation".

The paper is structured as follows. In Section \ref{S:canonical} we introduce a suitable family of standard domains and describe canonical mappings onto these domains in terms of the Green function and associated functions. In Section \ref{S:cr} we introduce an appropriate time parameter for the Jordan arcs in a multiply connected domain, namely the conformal radius, and establish a variational formula for ``increments" of the conformal radius under perturbations of the domain. In Section \ref{S:rkl} we establish what we call the radial Komatu-Loewner equation, which generalizes the radial Loewner equation to multiply connected domains. In Section \ref{S:mm} we obtain the corresponding equation for the evolution of the moduli and then show that the vector field in the differential equation  satisfies a Lipschitz property, Theorem \ref{T:lipschitz}. This result is crucial as it allows us to reverse the construction: Start with a motion on the boundary, solve the equation for the moduli for the given boundary motion, then solve the radial Komatu-Loewner equation to obtain a growing family of compacts.  In Section \ref{S:rSLE} we use the  correspondence between growing ``simple" curves and paths on moduli space to show that consistent conformally invariant random simple curves are given by diffusions on moduli space, and identify the diffusion up to a single drift term. Under an additional symmetry condition we then identify the diffusion up to a single parameter $\kappa$. In Section \ref{S:locality} we study the growing random compact sets for this diffusion and show that they satisfy the locality property if and only if $\kappa=6$.

As this paper neared its completion we became aware of the thesis of Dapeng Zhan, which contains another version of the radial Komatu-Loewner equation, and discusses a class of random Loewner chains on Riemann surfaces, \cite{zhan:2004}.  The results we present here were previously announced in \cite{bauer.friedrich:2004}.

The authors would like to acknowledge the hospitality of the Institut f\"ur Stochastik at the University of Freiburg and of the Isaac Newton Institute where part of this work was done.
R.F. would like to thank the Institute for Advanced Study for its hospitality and support. 

\section{Preliminaries}\label{S:prel}
When we generalize a problem in physics, usually additional degrees of freedom show up. In the case of the object one would like to call ``general SLE" one faces similar questions.

First, as we increase the connectivity of the domain, or look at arbitrary bordered surfaces, the possible simple paths of interest can {\it do more}, 
such that the notion of ``radial" and ``chordal" are not the only natural ones.

Let us illustrate this with the example of the annulus. A simple curve, starting at the outer boundary could either terminate at an interior point, connect the same boundary component or end at the inner boundary. Therefore we have to make a choice at the very beginning by restricting the possible behavior of the curves. 

However, there is a second component. As we are ultimately interested in random paths on surfaces, we have to define a probability measure on them. Since this is not given by geometric considerations, additional data is needed.

The natural point of view, at least from the standpoint of physics, is to see the curves as some characteristic manifestation of fields, defined on a (bordered) surface, which in the path-integral formulation, would be random, or quantum fields. A familiar geometric example of curves obtained via fields, are geodesics which are ``manifestations" of the symmetric non-degenerate two-tensor field, i.e. the metric. The field perspective naturally connects with CFT, at least if we are interested in conformally covariant properties, and opens up the door for the construction of a measure.  But to obtain from a field theory the desired class (in the geometric sense) of random curves, we have to choose boundary conditions. For a continuum theory the natural ones are either Neumann or Dirichlet and combinations of the two.

Again, in the familiar picture of statistical mechanics models on planer domains of some connectivity, such (random) curves would correspond  to domain walls, connecting points where the boundary conditions change discontinuously. 
This is a ``static" way of constructing measures on paths on surfaces, i.e. via the path-integral formalism (partition function).

However, by the use of the Loewner mapping we can transfer the problem from the static point of view to a dynamic one by characterizing the arising curves via their driving function.

Such a fundamental  approach was first introduced by Oded Schramm in~\cite{schramm:2000}  in order to describe the scaling limit of two specific models, namely the loop-erased random walk and the uniform spanning tree.

The natural assumptions for the above implied Brownian motion as driving function.

But a closer look reveals, that the most natural driving function for models arising from say statistical mechanics, is to assume a Markov process, which allows for an additional drift.

The crucial and  new insight is, that the link between the set of boundary conditions, which generate the desired class of random curves, and the corresponding Markov process is given by representation theoretic arguments. In the language of CFT, at the position where the boundary conditions change, we can insert boundary condition changing operators, i.e. local fields $|h\rangle$ of some real weight $h$.

As Hadamard's principle of boundary variation has a natural description in terms of operators, which one can show to satisfy the commutation relations of the Virasoro algebra, the singular deformation of the boundary is given by the operators $L_{-2}$ and $L_{-1}$ which correspond to ``cutting" the surface and ``moving" the point, where the boundary conditions change, i.e. the ``marked point". 

The condition on a vanishing drift as well as the way to calculate it at all, is exemplified in the following expression
\begin{displaymath}
\left(\frac{\kappa}{2}L^2_{-1}-2 L_{-2}\right)|h\rangle=0~.
\end{displaymath}

As we already mentioned in the introduction, another novel feature of the general Loewner evolution is the fact, that the driving function now lives on a higher-dimensional parameter space, i.e. a moduli space.

Finally we remark, that for a multiply connected domain, we still could encode the random curves by a random driving function on a one-dimensional space, but then the resulting stochastic process would not be Markov anymore, i.e. the path would have to remember its complete past, and the stochastic process would be given by a Girsanov transform.


\section{Standard domains and canonical mapping}\label{S:canonical}

Denote by $D$ a domain in the complex $z$-plane bounded by a finite
number of proper continua $C_j, j=1,\dots, n$, and let $w$ be a point in $D$. If $D$ is simply connected, then there is a unique conformal map $\Phi(z)=\Phi(z,w)$ from $D$ onto the unit disk such that 
\begin{equation}\label{E:normal}
	\Phi(w,w)=0,\quad 
	\frac{\partial \Phi(z,w)}{\partial z}|_{z=w}>0,
\end{equation}
If $n>1$, then there is a unique conformal map from $D$ onto the unit disk with $n-1$ disjoint concentric circular slits which maps $C_n$ to $\partial\mathbb D$ and satisfies \eqref{E:normal}. We call the unit disk with a finite number of concentric circular slits a {\em standard domain\ } and the normalized mapping $\Phi(z,w)$ the {\em canonical mapping}. We now recall the construction of the canonical mapping. 

For $z,w\in D$
denote $G_D(z,w)=G(z,w)$ the Green function of $D$ with pole at $w$.
$G(z,w)$ is a harmonic function of $z$ throughout $D$ except at $w$,
where $G(z,w)+\ln|z-w|$ is harmonic. Further, if $z$ converges to a
boundary continuum $C_j$, then $G(z,w)$ converges to zero. These
properties determine $G$.

If $D$ is simply connected, denote $H(z,w)$ a harmonic conjugate of
$G(z,w)$ with respect to $z$. $z\mapsto H(z,w)$ is multiple-valued (as $z$ describes a small circle around $w$, $H(z,w)$ changes by $2\pi$) and contains an arbitrary constant. If we choose the constant such that for one branch, and $x$ real, $\lim_{x\to0}H(x+w,w)=0$, then the canonical mapping is given by 
\begin{equation}\label{E:canonical1}
    \Phi(z,w)=\exp(-G(z,w)-i H(z,w)).
\end{equation}

If $n>1$, then $G$ has also periods with respect to circuits around $C_j$, $j=1,\dots, n-1$, 
\begin{equation}\label{E:harmonic measure}
	2\pi i\omega_j(z)=-i\int_{C_j}\frac{\partial G(z,w)}{\partial n_w}\ 
	ds_z.
\end{equation}
Here $\partial/\partial n$ denotes the derivative in the direction of the outward pointing normal and $ds$ denotes arc-length measure. For the purpose of this definition we have assumed that $C_j$ is piecewise smooth. By the Riemann mapping theorem and the conformal invariance of the Green function this represents no loss of generality. $\omega_j(z)$ is harmonic in $D$ with boundary values $1$ on $C_j$, $0$ on $C_k, k\neq j$, and is called the harmonic measure of $C_j$ in $z$ with respect to $D$. It is also the probability that planar Brownian motion started in $z$ exits $D$ through $C_j$. Denote $F(z,w)$ an analytic function in $z$ with real part $G(z,w)$. Then $\omega_j(z)$ is the real part of the analytic function
\begin{equation}\label{E:analytic hm}
	R_j(z)=-\frac{1}{2\pi}\int_{C_j}\frac{\partial F(z,w)}{\partial n_w}\ 	ds_w.
\end{equation}
$R_j(z)$ is regular in $D$ and possesses periods with respect to circuits around $C_k$ given by
\begin{equation}\label{E:periods}
	2\pi i P_{kj}=i\int_{C_k}\frac{\partial \omega_j(z)}{\partial n_z}\ 	ds_z.
\end{equation}
It is well known that the period matrix ${\bf P}=[P_{kj}]_{k,j=1}^{n-1}$ is  symmetric and positive-definite. Define the function
\begin{equation}\label{E:log map}
	F(z,w)+{\bf R}(z)^T {\bf P}^{-1}\boldsymbol\omega(w),
\end{equation}
where 
\[
	{\bf R}(z)^T=(R_1(z),\dots, R_{n-1}(z)),\quad 	\boldsymbol\omega(w)^T=(\omega_1(w),\dots,\omega_{n-1}(w)).
\]
This function has vanishing periods about $C_j,j=1,\dots,n-1$, and period $-2\pi$ about $C_n$. It is now easy to see that, after adding an appropriate imaginary constant, the canonical mapping for $D$ is given by
\begin{equation}\label{E:canonical n}
	\Phi(z,w)=\exp( -F(z,w)-{\bf R}(z)^T {\bf P}^{-1}\boldsymbol\omega(w)),
\end{equation}
see \cite{schiffer:1946} for details.

We can use $\Phi(z,w)$ to produce mappings onto other families of standard domains, see \cite{courant:1950}. For example, if $w=u+iv$, then 
\[
	z\mapsto -\frac{\partial}{\partial u}\ln\Phi(z,w)
\]
maps $D$ onto the whole plane slit along $n$ segments parallel to the imaginary axis. The point $z=w$ corresponds in this map to the point at infinity. Similarly,
\[
	z\mapsto -\frac{\partial}{\partial v}\ln\Phi(z,w)
\]
maps $D$ conformally onto the plane slit along segments parallel to the real axis. By letting $w$ approach one of the boundary components, say $C_n$, we find that 
\begin{equation}\label{E:standard2}
	z\mapsto\Psi(z,w):=\frac{\partial}{\partial n_w}\ln\Phi(z,w)
	=-\frac{\partial F(z,w)}{\partial n_w}
	-{\bf R}(z)^T {\bf P}^{-1}
	\frac{\partial\boldsymbol\omega(w)}{\partial n}
\end{equation}
is a conformal map from $D$ onto the right half-plane $\Re(\zeta)>0$ slit along $n-1$ segments parallel to the imaginary axis. If the boundary components $C_j$ are Jordan arcs, then $\Psi$ extends continuously to the boundary.  Note that if $D$ is the unit disk and $w=e^{i\varphi}$, then $\Psi$ is given by
\[
	z\mapsto\frac{e^{i\varphi}+z}{e^{i\varphi}-z}.
\]   


\section{Domain constant and conformal radius}\label{S:cr}

We define the {\em domain constant\ } $d_D(w)$ by
\begin{equation}  \label{E:domain constant}
	d_D(w)=-\lim_{z\to w}\left(G(z,w)+\ln|z-w|\right),
\end{equation}
and the conformal radius $r_D(w)$ by
\begin{equation}\label{E:conformal radius}
	r_D(w)=\ln \left(\frac{\partial\Phi(z,w)}{\partial z}|_{z=w}\right).
\end{equation}
If $D$ is simply connected then we find from \eqref{E:canonical1} that $d_D(w)=r_D(w)$. If $n>1$, then it follows from \eqref{E:canonical n} that 
\begin{equation}\label{E:conformal radius 2}
	r_D(w)=-\lim_{z\to w}\left(\gamma(z,w)+\ln|z-w|\right),
\end{equation}
where 
\begin{equation}\label{E:gamma}
	\gamma(z,w):=\Re(\ln\Phi(z,w))=G(z,w)+\boldsymbol\omega(w)^T{\bf P^{-1}}\boldsymbol\omega(z).
\end{equation}
We note that if $D$ is a standard domain, then $\gamma(z,0)=-\ln|z|$ and thus  
\begin{equation}\label{E:green-gamma}
	 	G(z,0)=-\ln|z|-\boldsymbol\omega(0)^T{\bf P^{-1}}\boldsymbol\omega(z).
\end{equation}
It follows from \eqref{E:conformal radius 2}, \eqref{E:domain constant}, and \eqref{E:gamma} that 
\begin{equation}\label{E:dcn}
	d_D(w)=r_D(w)+\omega(w)^T {\bf P}^{-1}\omega(w).
\end{equation}

We call a closed and simply connected set $A$ a {\em hull in} $D$
if $\overline{A\cap D}=A$ and $\mathbb D\backslash A$ has the same connectivity as $D$. Suppose that $E$ is a standard domain and consider two disjoint hulls in $E$, say
$A$ and $B$, that do not contain $0$. Denote $\Phi_A(z)=\Phi_A(z,0)$ the
canonical mapping from $E\backslash A$ onto the standard domain $E^*:=\Phi_A(E\backslash A)$ that fixes zero. Let $D=E\backslash B$, $D^*=E^*\backslash \Phi_A(B)$, and set $\Gamma=\partial D$,
$\Gamma^*=\partial D^*$. We want to compare the ``domain constant increments"
$d_{D^*}(0)-d_{D}(0)$ and $d_{E^*}(0)-d_E(0)$. Assume that $B$ is connected, intersects the unit circle $\partial \mathbb D$, and has a piecewise smooth
boundary. Since $G_{D^*}(z,w)-G_{E^*}(z,w)$ is
harmonic throughout $D^*$, the Poisson formula gives
\begin{equation}\label{E:hadamard1}
    G_{D^*}(z,w)-G_{E^*}(z,w) 				  	=\frac{1}{2\pi}\int_{\Gamma^*}(G_{D^*}(\eta,w)-	G_{E^*}(\eta,w)) 
    \frac{\partial G_{D^*}(\eta,z)}{\partial n_{\eta}}\ ds_{\eta},
\end{equation}
where $\partial/\partial n$ denotes the derivative along the outward
pointing normal and $ds$ integration relative to arc length. Note that we integrate along both sides of the slits. Since $G_D^*(\eta,w)=0$ for $\eta\in\Gamma^*$, and $G_{E^*}(\eta,w)=0$ for $\eta\in\Gamma^*\backslash\Phi_A(B)$, we get from\eqref{E:hadamard1} and \eqref{E:domain constant}
\begin{equation}\label{E:inc1}
	d_{D^*}(0)-d_{E^*}(0)=\frac{1}{2\pi}\int_{\Phi_A(\partial B)\cap E^*}
	G_{E^*}(\eta,0)\frac{\partial G_{D^*}(\eta,0)}{\partial n_{\eta}}\ 	ds_{\eta}.
\end{equation}
Similarly,
\begin{equation}\label{E:inc2}
	d_{D}(0)-d_{E}(0)=\frac{1}{2\pi}\int_{\partial B\cap E}
	G_{E}(\eta,0)\frac{\partial G_{D}(\eta,0)}{\partial n_{\eta}}\ 	ds_{\eta}.
\end{equation}
Let now $\{B_{\epsilon},\epsilon>0\}$ be a family of hulls as $B$ above,
such that $B_{\epsilon}\supset B_{\epsilon'}$ if
$\epsilon>\epsilon'$, and $diam(B_{\epsilon})=O(\epsilon)$. Let
$\xi=\bigcap_{\epsilon>0}B_{\epsilon}$. We note that necessarily
$\xi\in\partial E^*$. In the following we suppress the subscript $\epsilon$ when there is no risk of confusion. For $\eta\in\partial B_{\epsilon}\cap E$ we have by Hadamard's formula
\begin{align}\label{E:hadamard}
	\frac{\partial G_{D}(\eta,0)}{\partial n_{\eta}}
	&=\frac{\partial G_{E}(\xi,0)}{\partial n_{\eta}}+o(1),\notag\\
	\frac{\partial G_{D^*}(\Phi_A(\eta),0)}{\partial n}
	&=\frac{\partial G_{E^*}(\Phi_A(\xi),0)}{\partial n_{\eta}}+o(1),
\end{align}
and by Taylor's formula
\begin{align}\label{E:taylor}
	G_{E}(\eta,0)
	&=\left(-\frac{\partial G_{E}(\xi,0)}{\partial n_{\eta}}+O(\epsilon)\right)|\eta-\xi|,
	\notag\\
	G_{E^*}(\Phi_A(\eta),0)
	&=\left(-\frac{\partial G_{E^*}(\Phi_A(\xi),0)}{\partial 	n_{\eta}}+O(\epsilon)\right)|\Phi_A(\eta)-\Phi_A(\xi)|,
\end{align}
 and also
\begin{equation}\label{E:taylor2}
	\int_{\Phi_A(\partial B\cap E)}|\eta-\Phi_A(\xi)|\ 	ds_{\eta}=|\Phi_A'(\xi)|^2\int_{\partial B\cap E}|\eta-\xi|\ 	ds_{\eta}+o(\epsilon^2).
\end{equation}
Thus, from \eqref{E:inc1}, \eqref{E:inc2}, \eqref{E:hadamard}, \eqref{E:taylor}, and \eqref{E:taylor2}, we get finally
\begin{equation}\label{E:change of d}
	\lim_{\epsilon\to0}\frac{d_{D^*}(0)-d_{E^*}(0)}{d_D(0)-d_E(0)}
	=\left(\frac{\frac{\partial G_{E^*}(\Phi_A(\xi),0)}{\partial 	n_{\eta}}}{\frac{\partial G_{E}(\xi,0)}{\partial n_{\eta}}}\right)^2|
	\Phi_A'(\xi)|^2.
\end{equation}

If the hulls $B_{\epsilon}$ are not smooth, then, given any
$\delta>0$, we may approximate them by smooth hulls
$B_{\epsilon,\delta}$, such that
\[
    1-\delta\le
    \frac{d_{D_{\delta}}(0)-d_{E_{\delta}}(0)}{d_D(0)-d_E(0)},\frac{ 	d_{D^*_{\delta}}(0)-d_{E^*_{\delta}}(0)}{d_{D^*}(0)-d_{E^*}(0)}\le
    1+\delta.
\]
Hence equation \eqref{E:change of d} applies also in the non-smooth
case.

Next, we wish to compare the ``domain constant increment" with the ``conformal radius increment." By \eqref{E:dcn}
\begin{equation}\label{E:d over r}
	\frac{d_D(0)-d_E(0)}{r_D(0)-r_E(0)}=1+\frac{
	\boldsymbol\omega_D(0)^T{\bf P}_D^{-1}\boldsymbol\omega_D(0)
	-\boldsymbol\omega_E(0)^T{\bf P}_E^{-1}\boldsymbol\omega_E(0)}{r_D(0)-r_E(0)}.
\end{equation}
We consider the domain $D$ as a variation of the domain $E$, i.e. we continue to consider the family $\{B_{\epsilon},\epsilon>0\}$ from above. From Hadamard's formula, \cite[(45)]{schiffer:1946} and \cite[Chap. I, Sec. 11]{nehari:1952}, we have
\begin{align}
	\lim_{\epsilon\to0}\frac{r_D(0)-r_E(0)}{\epsilon}&=1\label{E:hcr}\\
	\lim_{\epsilon\to0}\frac{
	\boldsymbol\omega_D(0)-\boldsymbol\omega_E(0)}{\epsilon}&=
	-\frac{
	\partial\boldsymbol\omega_E(\xi)}{\partial n}\cdot\frac{
	\partial G_E(\xi,0)}{\partial n_1}\\
	\lim_{\epsilon\to0}\frac{{\bf P}_D-{\bf P}_E}{\epsilon}&=\frac{
	\partial\boldsymbol\omega(\xi)}{\partial n}\cdot\frac{
	\partial\boldsymbol\omega(\xi)}{\partial n}^T.
\end{align} 
By \eqref{E:green-gamma}, it follows that
\begin{equation}\label{E:hadamard hm}
	\lim_{\epsilon\to0}\frac{
	\boldsymbol\omega_D(0)-\boldsymbol\omega_E(0)}{\epsilon}
	=\frac{
	\partial\boldsymbol\omega_E(\xi)}{\partial n}+\frac{
	\partial\boldsymbol\omega_E(\xi)}{\partial n}\boldsymbol\omega_E(0)^T{\bf P}_E^{-1}\frac{
	\partial\boldsymbol\omega_E(\xi)}{\partial n}.
\end{equation}
We also have
\begin{equation}\label{E:hadamard P-1}
	\lim_{\epsilon\to0}\frac{{\bf P}_D^{-1}-{\bf P}_E^{-1}}{\epsilon}=-{\bf P}_E^{-1}\frac{
	\partial\boldsymbol\omega(\xi)}{\partial n}\cdot\frac{
	\partial\boldsymbol\omega(\xi)}{\partial n}^T {\bf P}_E^{-1}.
\end{equation}
Combining \eqref{E:d over r}, \eqref{E:hcr}, \eqref{E:hadamard hm}, and \eqref{E:hadamard P-1} we get
\begin{align}\label{E:d over r2}
	\lim_{\epsilon\to0}\frac{d_D(0)-d_E(0)}{r_D(0)-r_E(0)}
	&=\left(1+\boldsymbol\omega_E(0)^T{\bf P}_E^{-1}\frac{
	\partial\boldsymbol\omega_E(\xi)}{\partial n}\right)^2\notag\\
	&=\left(\frac{\partial G_E(\xi,0)}{\partial n}\right)^2.
\end{align}
From \eqref{E:change of d} and \eqref{E:d over r2} it now follows that 
\begin{equation}\label{E:change of r}
	\lim_{\epsilon\to0}\frac{r_{D^*}(0)-r_{E^*}(0)}{r_D(0)-r_E(0)}
	=|\Phi_A'(\xi)|^2.
\end{equation}


\section{Radial Komatu-Loewner equation}\label{S:rkl}

In this section we derive what we call the {\em radial Komatu-Loewner equation.} It describes the evolution of slit mappings when the slit grows from a point on the boundary of a multiply connected domain to a point in the interior. It shows in particular, that each growing slit induces a continuous motion on one boundary component of the domain, namely the image of the tip of the slit under the canonical mapping. The radial Loewner equation, of course, is the special case when the domain is simply connected. 

In \cite{komatu:1950}, Komatu derives a differential equation that is satisfied by the canonical mappings for the growing slit in the case where the slit grows from one boundary component to another boundary component. The radial Komatu-Loewner equation is similar to the equation Komatu derives in \cite{komatu:1950} and our proof proceeds along the lines of the proof he gave. However, the equation given in \cite{komatu:1950} involves the derivatives of certain moduli with respect to the slit-parameter and it is only stated that these derivatives exist, but  their explicit form is not given. Since we wish to solve the radial Komatu-Loewner equation for a given input, we need the explicit form of the equation. Using the fact that the real part of a single-valued analytic function is orthogonal to the real parts of Abelian differentials of the first kind, we derive the explicit equation directly. 

For the purposes we have in mind---generating random slits---the main difference between the simply and the multiply connected case is, that in the multiply connected case a growing slit corresponds to a motion on the boundary coupled to a motion on the moduli space. In the simply connected case the moduli space reduces to a single point.  

Let $D$ be a standard domain, and $\gamma:[0,t_{\gamma}]\to\overline{D}$ a Jordan arc such that $\gamma(0)\in S^1$, and $\gamma(0,t_{\gamma}]\subset D\backslash\{0\}$. Let $g_t$ be the canonical mapping from $D\backslash\gamma[0,t]$ which leaves $0$ fixed, and denote $D_t$ the standard domain $g_t(D\backslash\gamma[0,t])$. By \eqref{E:hcr}, the map $t\mapsto g_t'(0)$ is a strictly increasing continuous function, see \cite{schiffer:1946} for details, and we may assume that the parameter $t$ is chosen so that $g_t'(0)=e^t$. We wish to find a differential equation for the family $\{g_t:t\in[0,t_{\gamma}]\}$. 

Denote $C_j(t), j=1,\dots,n$, the boundary components of $D_t$. We always have $C_n(t)=S^1$. For $j=1,\dots, n-1$, let $m_j(t)$ be the radial distance of the circular slit $C_j(t)$ from the origin. Denote $\xi(t)$ the starting point on $S^1$ of the Jordan arc $g_t(\gamma[t,t_{\gamma}])$ in $D_t$, i.e. $g_t(\gamma_t)$. For $0<t^*<t<t_{\gamma}$, set 
\[
	g_{t,t^*}=g_{t^*}\circ g_t^{-1}.
\]  
Then $g_{t,t^*}$ is a conformal map from $D_t$ onto $D_{t^*}\backslash g_{t^*}(\gamma[t^*,t])$. The point $\xi(t^*)=g_{t^*}(\gamma_{t^*})$ corresponds to two prime ends in  $D_{t^*}\backslash g_{t^*}(\gamma[t^*,t])$. Denote $\exp(i\beta_0(t,t^*))$ and $\exp(i\beta_1(t,t^*))$, with $\beta_0(t,t^*)<\beta_1(t,t^*)$, the pre-images of these prime ends under $g_{t,t^*}$, i.e.
\[
	g_{t,t^*}(\exp(i\beta_0(t,t^*)))=g_{t,t^*}(\exp(i\beta_1(t,t^*)))
	=g_{t^*}(\gamma_{t^*}).
\]
Then, if $|z|=1$ and $\beta_1(t,t^*)\le\arg z\le\beta_0(t,t^*)+2\pi$, 
\[
	|g_{t,t^*}(z)|=1.
\]

Consider the function 
\[
	z\mapsto\ln\frac{g_{t,t^*}(z)}{z}.
\]
Since $g_{t,t^*}(0)=0$, it is regular throughout $D_t$. In particular, 
\begin{equation}\label{E:zerovalue}
	\lim_{z\to0}\frac{g_{t,t^*}(z)}{z}=e^{t^*-t}.
\end{equation}
Furthermore, as $z$ describes a small circle around $0$, $\ln z$ changes to $\ln z+2\pi i$ and $\ln g_{t,t^*}(z)$ changes to $\ln g_{t,t^*}(z)+2\pi i$. Thus $\ln (g_{t,t^*}(z)/z)$ is regular and one-valued throughout $D_t$. As $\ln|g_{t,t^*}(z)/z|$ is regular and harmonic in $D_t$, Poisson's formula gives
\begin{equation}\label{E:poisson}
	\ln\left|\frac{g_{t,t^*}(z)}{z}\right|=-\frac{1}{2\pi}\int_{\partial D_t}
	\ln\left|\frac{g_{t,t^*}(\zeta)}{\zeta}\right|\frac{\partial G(\zeta,z;t)}{\partial n_1}\ ds,
\end{equation}
where $G(\zeta,z;t)$ is the Green function for $D_t$ with pole at $z$. We write $u$ for the harmonic function of $z$ on the left-hand side of \eqref{E:poisson} and denote by $v$ its harmonic conjugate. Since $u$ is the real part of a one-valued analytic function, the period $P_j(v)$ of $v$ with respect to a circuit about $C_j(t)$ has to vanish for each $j=1,\dots, n-1$. Since $\omega_j(z;t)=\delta_{j,k}$ for  $z\in C_k(t)$, we have
\begin{align}\label{E:orthogonal}
	0=P_j(v)&=\int_{C_j(t)} dv
	=\int_{C_j(t)}\frac{\partial v}{\partial s}\ ds
	=\int_{C_j(t)}\frac{\partial u}{\partial n}\ ds\notag\\
	&=\int_{\partial D_t}\omega_j(\zeta;t)\frac{\partial u}{\partial n}\ ds 	=\int_{\partial D_t}u\frac{\partial\omega_j(\zeta;t)}{\partial n}\ ds. 
\end{align} 
Combining \eqref{E:orthogonal} and \eqref{E:poisson}, we find  
\begin{align}\label{E:poisson2}
	\ln&\left|\frac{g_{t,t^*}(z)}{z}\right|\notag\\
	&=-\frac{1}{2\pi}\int_{\partial D_t}
	\ln\left|\frac{g_{t,t^*}(\zeta)}{\zeta}\right|
	\left(\frac{\partial G(\zeta,z;t)}{\partial n_1}
	+\boldsymbol{\omega}(z;t)^T\boldsymbol{P}_t^{-1}\frac{\partial\boldsymbol{\omega}(\zeta;t)}{\partial n}\right) ds.
\end{align}
It follows from Section \ref{S:canonical} that 
\[
	z\mapsto \frac{\partial G(\zeta,z;t)}{\partial n_1}
	+\boldsymbol{\omega}(z;t)^T\boldsymbol{P}_t^{-1}\frac{\partial\boldsymbol{\omega}(\zeta;t)}{\partial n}
\]
has a single-valued harmonic conjugate and so
\begin{align}\label{E:complex poisson}
	&\ln\frac{g_{t,t^*}(z)}{z}\notag\\
	&=-\frac{1}{2\pi}\int_{\partial D_t}
	\ln\left|\frac{g_{t,t^*}(\zeta)}{\zeta}\right|\left(\frac{\partial F(\zeta,z;t)}{\partial n_1}+\boldsymbol{R}(z;t)^T\boldsymbol{P}_t^{-1}\frac{\partial\boldsymbol{\omega}(\zeta;t)}{\partial n}\right)\ ds\notag\\
	&\qquad +ic,
\end{align}
where $z\mapsto F(\zeta,z;t)$ is a multiple-valued analytic function with real part $G(\zeta,z;t)$, and where $c$ is a real constant. Note that by \eqref{E:analytic hm} and \eqref{E:periods}, for every $j=1,\dots,n-1$,
\begin{align}
	-\frac{1}{2\pi}&\int_{C_j(t)}
	\left(\frac{\partial F(\zeta,z;t)}{\partial n_1}
	+\boldsymbol{R}(z;t)^T\boldsymbol{P}_t^{-1}
	\frac{\partial\boldsymbol{\omega}(\zeta;t)}{\partial n}\right)\ ds	\notag\\
	&=R_j(z;t)-\boldsymbol{R}(z;t)^T\boldsymbol{P}_t^{-1}(\boldsymbol{P}_t)_{.,j}=0.
\end{align}
Thus, since  $\ln|g_{t,t^*}(\zeta)/\zeta|$ is constant on each $C_j(t)$, $j=1,\dots, n-1$, and identically zero on $\{|z|=1,\beta_1(t,t^*)\le\arg z\le\beta_0(t,t^*)+2\pi\}$, it follows from \eqref{E:complex poisson} that 
\begin{align}\label{E:complex poisson2}
	&\ln\frac{g_{t,t^*}(z)}{z}\notag\\
	&=-\frac{1}{2\pi}\int_{\beta_0(t,t^*)}^{\beta_1(t,t^*)}
	\ln\left|\frac{g_{t,t^*}(\zeta)}{\zeta}\right|\left(\frac{\partial F(\zeta,z;t)}{\partial n_1}+\boldsymbol{R}(z;t)^T\boldsymbol{P}_t^{-1}\frac{\partial\boldsymbol{\omega}(\zeta;t)}{\partial n}\right)\ ds \notag\\
	&\qquad+ic.
\end{align}
We now show that $c=0$. To that end, note first that by Cauchy's integral formula,
\begin{equation}\label{E:cauchy}
	t^*-t=\frac{1}{2\pi i}\int_{\partial D_t}
	\ln\left(\frac{g_{t,t^*}(\zeta)}{\zeta}\right)\ 	\frac{d\zeta}{\zeta}.
\end{equation}
In particular, the right-hand side of \eqref{E:cauchy} is real.  Since all boundary components are concentric circular slits, $d\zeta/\zeta$ is purely imaginary along $\partial D_t$, i.e.
\[
	\frac{d\zeta}{\zeta}=i\ d\arg\zeta,\quad\zeta\in\partial D_t.
\]
Hence
\begin{align}
	t^*-t&=\frac{1}{2\pi}\int_{\partial D_t}
	\ln\left|\frac{g_{t,t^*}(\zeta)}{\zeta}\right|\ d\arg\zeta\notag\\
	&=\frac{1}{2\pi}\int_{\beta_0(t,t^*)}^{\beta_1(t,t^*)}
	\ln\left|g_{t,t^*}(e^{i\varphi})\right|\ d\varphi 	+\frac{1}{2\pi}\sum_{j=1}^{n-1}
	\int_{C_j(t)}\ln\frac{m_j(t^*)}{m_j(t)}\ d\arg\zeta.
\end{align}
Since the two ``sides" of $C_j(t)$ make opposite contributions,
\[
	\int_{C_j(t)}d\arg\zeta=0,\quad j=1,\dots, n-1,
\]
and we finally get
\begin{equation}\label{E:increment}
	t^*-t=\frac{1}{2\pi}\int_{\beta_0(t,t^*)}^{\beta_1(t,t^*)}
	\ln\left|g_{t,t^*}(e^{i\varphi})\right|\ d\varphi.
\end{equation}
Note next that, by \eqref{E:canonical n},
\[
	F(\zeta,0;t)+\boldsymbol{R}(\zeta;t)^T\boldsymbol{P}_t^{-1}
	\boldsymbol{\omega}(0;t)=-\ln \zeta.
\]
Further, since we are free to add a constant to the harmonic conjugate, we may assume that $\Im(\boldsymbol{R}(0;t))=\boldsymbol{0}$, i.e that $\boldsymbol{\omega}(0;t)=\boldsymbol{R}(0;t)$. Also, since $\boldsymbol\omega(\zeta;t)$ is constant along the boundary, 
\[
	\frac{\partial\Im(\boldsymbol{R}(\zeta;t)}{\partial n}=\frac{\partial\boldsymbol{\omega}(\zeta;t)}{\partial s}=0.
\]
Thus, using $\boldsymbol P=\boldsymbol P^T$, we get
\begin{align}\label{E:what}
	-1&=-\frac{\partial\ln\zeta}{\partial n} \notag\\
	&=\frac{\partial F(\zeta,0;t)}{\partial n_1}+\boldsymbol{R}(0;t)^T\boldsymbol{P}_t^{-1}\frac{\partial\boldsymbol{\omega}(\zeta;t)}{\partial n}.
\end{align}
Hence, if we evaluate \eqref{E:complex poisson2} at $z=0$ and use \eqref{E:zerovalue} on the left, and \eqref{E:increment}, \eqref{E:what} on the right, then it follows that $c=0$.

Letting $z=g_t(w)$ we now get
\begin{align}\label{E:difference}
	\ln&\frac{g_{t^*}(w)}{g_t(w)}\notag\\
	&=-\frac{1}{2\pi}\int_{\beta_0(t,t^*)}^{\beta_1(t,t^*)}
	\ln|g_{t,t^*}(e^{i\varphi)})|\notag\\
	&\qquad\qquad\times\left(\frac{\partial F(e^{i\varphi},z;t)}{\partial n_1}
	+\boldsymbol{R}(z;t)^T\boldsymbol{P}_t^{-1}
	\frac{\partial\boldsymbol{\omega}(e^{i\varphi};t)}{\partial n}\right)\ ds.
\end{align}

We now wish to let $t^*\nearrow t$ in \eqref{E:difference}. Note first that, for $\varphi\in[0,2\pi]$, $\varphi\mapsto\ln|g_{t,t^*}(e^{i\varphi})|$ is continuous and non-positive and that also
\[
	\varphi\mapsto A(\varphi):=
	\frac{\partial F(e^{i\varphi},z;t)}{\partial n_1}
	+\boldsymbol{R}(z;t)^T\boldsymbol{P}_t^{-1}
	\frac{\partial\boldsymbol{\omega}(e^{i\varphi};t)}{\partial n}
\]
is continuous. Thus it follows from the mean-value theorem of integration that 
\begin{align}
	\frac{1}{2\pi(t^*-t)}&\int_{\beta_0(t,t^*)}^{\beta_1(t,t^*)}
	\ln\left|g_{t,t^*}(e^{i\varphi})\right| A(\varphi)\ d\varphi\notag\\
	&=\frac{\Re(A(\varphi'))+i\Im(A(\varphi''))}{2\pi(t^*-t)}
	\int_{\beta_0(t,t^*)}^{\beta_1(t,t^*)}
	\ln\left|g_{t,t^*}(e^{i\varphi})\right|\ d\varphi\notag\\
	&=\Re(A(\varphi'))+i\Im(A(\varphi'')),
\end{align}
for some $\varphi',\varphi''\in[\beta_0(t,t^*),\beta_1(t,t^*)]$. Hence
\begin{align}
	\lim_{t^*\nearrow t}&\frac{\ln g_{t^*}(w)-\ln g_t(w)}{t^*-t}\notag\\
	&=-\frac{\partial F(\xi(t),z;t)}{\partial n_1}
	-\boldsymbol{R}(z;t)^T\boldsymbol{P}_t^{-1}
	\frac{\partial\boldsymbol{\omega}(\xi(t);t)}{\partial n}.
\end{align}
By the same argument we may let $t\searrow t^*$. On the right-hand side above we then only need to change $t$ to $t^*$ and introduce an overall minus sign. Thus we have established the following

\begin{theorem}[Radial Komatu-Loewner equation]
	If $\gamma$ is a Jordan arc in a standard domain $D$ starting on $S^1$ with the parametrization from above, and if $g_t$ is the canonical map for $D\backslash\gamma[0,t]$, then, using the notation from above, the family $\{g_t:t\in[0,t_{\gamma}]\}$ satisfies the equation
\begin{equation}\label{E:RKL}
\partial_t\ln g_t(z)
	=-\frac{\partial F(\xi(t),g_t(z);t)}{\partial n_1}
	-\boldsymbol{R}(g_t(z);t)^T\boldsymbol{P}_t^{-1}
	\frac{\partial\boldsymbol{\omega}(\xi(t);t)}{\partial n},
\end{equation} 
with initial condition $g_0(z)=z$.
\end{theorem}


\section{Motion of moduli}\label{S:mm}

The right-hand side of the radial Komatu-Loewner equation, at time $t$, involves the Green function of the domain $D_t$, and also various functions derived from the Green function. Consequently, it does not make sense to ask for the solution of \eqref{E:RKL} for a given continuous curve $t\mapsto\xi(t)$, since the vector-field on the right-hand side of \eqref{E:RKL} is not specified by giving that information alone. To specify the Green function of $D_t$ we also need the {\em moduli} of the domain $D_t$. We will now consider what the appropriate moduli space is for our purposes and find a system of equations these moduli satisfy. Once this system is found, we can solve it for a given input $t\mapsto\xi(t)$, and then, in a second step, solve the radial Komatu-Loewner equation using $\xi$ and the moduli.

The geometric description of $D_t$ requires $3n-3$ real parameters, three for each (interior) slit, given, for example, by the radial distances of the slits to the origin, i.e $m_j(t)$, $j=1,\dots, n-1$, and the angles \[
	\theta_j(t)<\theta_j'(t),\quad j=1,\dots,n-1,
\]
determining the endpoints of the slit $C_j(t)$, $j=1\dots,n-1$. On the other hand, it is well known that two $n$-connected domains with non-degenerate boundary continua are conformally equivalent if $3n-6$ real parameters agree for $n>2$. If $n=2$ then there is only one real parameter describing the conformal class, and if $n=1$, then all such domains are conformally equivalent. In our context, we only allow conformal maps for which a certain interior point has a prescribed image and whose derivative at that point is positive---the canonical maps from Section \ref{S:canonical}. This corresponds to considering domains with one marked interior point and one marked boundary point. Indeed, given a domain $D$ with $n$ non-degenerate boundary continua, one interior point $z$ and one  boundary point (or, more generally, prime end) $\zeta$, there is a unique conformal map from $D$ onto a standard domain such that $z$ is mapped to $0$, and $\zeta$ to $1$. 

For the slits we wish to grow the marked points are the beginning ($t=0$) and end points $(t=\infty)$. It is now easy to see that the moduli space of $n$-connected planar domains with one interior and one boundary point marked is $3n-3$ dimensional for all $n\ge1$. We will take $\boldsymbol{m}(t)=(m_1(t),\dots,m_{n-1}(t))$, $\boldsymbol{\theta}(t)=(\theta_1(t),\dots,\theta_{n-1}(t))$, and $\boldsymbol{\theta}'(t)=(\theta_1'(t),\dots,\theta_{n-1}'(t))$ as the moduli of the domain $D_t$ and write $\boldsymbol{M}(t):=(\boldsymbol{m}(t),\boldsymbol{\theta}(t),\boldsymbol{\theta}'(t))$. For a standard domain the marked points are 0 and 1. For a point $\boldsymbol M$ in the moduli space we denote by $D=D({\boldsymbol M})$ the corresponding standard domain, while for an arbitrary $n$-connected domain $D$ we write $\boldsymbol M=M(D)$ for the corresponding point in the moduli space. 

Set
\begin{equation}\label{E:psi}
	\Psi_t(z,\zeta)
	=-\frac{\partial F(\zeta,z;t)}{\partial n_1}
	-\boldsymbol{R}(z;t)^T\boldsymbol{P}_t^{-1}
	\frac{\partial\boldsymbol{\omega}(\zeta;t)}{\partial n}.
\end{equation}
Then, by \eqref{E:standard2}, $\Psi_t$ maps $D_t$ conformally onto a slit right half-plane. By \eqref{E:what}, $\Psi_t$ is the unique such map with 
\begin{equation}\label{E:normalPsi}
	\Psi_t(0,\xi)=1.
\end{equation}
We may write the radial Komatu-Loewner equation \eqref{E:RKL} as
\[
	\partial_t\ln g_t(z)=\Psi_t(g_t(z),\xi(t)).
\]
By boundary correspondence, if $z\in C_j$, then $g_t(z)\in C_j(t)$ and \[
	\Re(\ln g_t(z))=\ln m_j(t).
\]
Thus, by considering the real part of the radial Komatu-Loewner equation,
\begin{equation}\label{E:realRKL}
	\partial_t\ln m_j(t)=\Re(\Psi_t(g_t(z),\xi(t))).
\end{equation}
Further, if 
\[
	z_j(t)=m_j(t) e^{i\theta_j(t)},\quad z_j'(t)=m_j(t)e^{i\theta_j'(t)}
\]
are the endpoints of the slit $C_j(t)$, then
\[
	z_j(t)=g_t(m_j(0)e^{i\eta_j(t)}),\quad z_j'(t)=g_t(m_j(0)e^{i\eta_j'(t)}),
\]
where $\theta_j(0)<\eta_j(t),\eta_j'(t)<\theta_j'(0)$. Indeed, the pre-images of the tips of $C_j(t)$, that is $m_j(0)e^{i\eta_j(t)}$ and $m_j(0)e^{i\eta_j'(t)}$, are the solutions to the equation
\[
	\frac{\partial}{\partial z}g_t(z)=0,
\]
on the set of prime-ends corresponding to $ C_j\backslash\{z_j(0),z_j'(0)\}$. A tip of $C_j(t)$ cannot be the image of a tip of $C_j$ because then the analytic function $\partial g_t/\partial z$ would not have the required number of zeroes, $2n-2$. 	
	
\begin{lemma}[Motion of moduli]\label{L:mm}
	The moduli 												\[														\boldsymbol{M}(t)=(\boldsymbol{m}(t),\boldsymbol{\theta}(t),
	\boldsymbol{\theta}'(t))
	\] satisfy the system of equations
	\begin{align}\label{E:mm}
		\partial_t\ln m_j(t)&=-\left(\boldsymbol{P}_t^{-1}
		\frac{\partial\boldsymbol{\omega}(\xi(t);t)}{\partial n} 		\right)_j,\notag\\
		\partial_t\theta_j(t)&=
		\Im\left(\Psi_t\left(m_j(t)e^{i\theta_j(t)},\xi(t)\right)\right),
		\notag\\
		\partial_t\theta_j'(t)&=
		\Im\left(\Psi_t\left(m_j(t)e^{i\theta_j'(t)},\xi(t)\right)\right),
	\end{align}
	for  $j=1,\dots,n-1$.
\end{lemma}

\begin{proof}
	We note that $\partial g_t/\partial z$ and $\partial^2 g_t/(\partial z)^2$ are analytic functions that extend analytically to the prime-ends corresponding to $C_1,\dots, C_{n-1}$ with the endpoints of the slits removed. By the implicit function theorem, 
\[
	t\mapsto m_j(0)e^{i\eta_j(t)}
\]
is differentiable with derivative 
\[
	DER_t:=\left[\frac{\partial^2 g_t}{(\partial z)^2}(m_j(0)e^{i\eta_j(t)})\right]^{-1}\frac{\partial^2 g_t}{\partial t\partial z}(m_j(0)e^{i\eta_j(t)}).
\]
By counting zeroes we find that 
\[
	\frac{\partial^2 g_t}{(\partial z)^2}(m_j(0)e^{i\eta_j(t)})\neq0
\]
and so $DER_t$ is finite. Hence
\begin{align}
	\partial_t\theta_j(t)&=\partial_t(\Im(g_t(m_j(0)e^{i\eta_j(t)}))\notag\\
	&=\Im\left(\Psi_t\left(m_j(t)e^{i\theta_j(t)},\xi(t)\right)\right)
	+\Im\left((\partial_z g_t)(z_j(t)) \times DER_t\right)\notag\\
	&=\Im\left(\Psi_t\left(m_j(t)e^{i\theta_j(t)},\xi(t)\right)\right).
\end{align}
In a similar way we obtain the derivative of $\theta_j'(t)$.
It remains to check that \eqref{E:realRKL} agrees with the 	first equation in \eqref{E:mm}. To this end we note that
	\[
		  \Re(\Psi_t(z,\zeta))
		  =-\frac{\partial G(\zeta,z;t)}{\partial n_1}
	-\boldsymbol{\omega}(z;t)^T\boldsymbol{P}_t^{-1}
	\frac{\partial\boldsymbol{\omega}(\zeta;t)}{\partial n}.
	\]
	From the boundary behavior of the Green function and the 	harmonic measures, it follows that for $z\in C_j(t)$
	\[
		\frac{\partial G(\zeta,z;t)}{\partial n_1}=0, 		\quad\text{and }\omega_k(z)=\delta_{jk}.
	\]
	The lemma follows.
\end{proof}


We now have our main existence statement.
	
\begin{theorem}\label{T:lipschitz}
	Given a continuous function $t\in[0,\infty)\mapsto\xi(t)\in S^1$ 	and the moduli $\boldsymbol{M}$ of a standard domain $D$, 	there exists a unique solution $\boldsymbol{M}(t)$ to the 	system \eqref{E:mm} on an interval $[0,t_{\xi})$ with 	$\boldsymbol{M}(0)=\boldsymbol{M}$, and where $t_{\xi}$ is 	characterized by 
	\[
		t_{\xi}=\inf\{\tau:\lim_{t\nearrow \tau}m_j(t)=1\text{ 
	for some }j\in\{1,\dots, n-1\}\}.
	\] 
	Further, if $D_t$ is the standard domain determined by 	$\boldsymbol{M}(t)$, and if $\Psi_t(z,\zeta)$ is the holomorphic 	vector field associated to $D_t$ by \eqref{E:psi}, then, for any 	$z\in D$, the equation
	\[	
		\partial_t \ln g_t^{D}(z)=\Psi_t(g_t^D(z),\xi(t)),\quad g_0^D(z)=z,
	\]
	has a unique solution on $[0,t_z)$, where 
	\[
		t_z=\sup\{t\le t_{\xi}:\inf_{s\in[0,t]}|g_s^D(z)-\xi(s)|>0\}.
	\]
	Finally, for $t<t_{\xi}$ set $K_t=\{z\in D:t_z\le t\}$. Then $g_t^D$ is the canonical conformal map from $D\backslash K_t$ onto $D_t$ which fixes zero and has positive derivative there.
\end{theorem}

\begin{proof}
For the existence of the solution to the moduli equations \eqref{E:mm} on $[0,t_\xi)$ we will show that the vector field in \eqref{E:mm} is Lipschitz as a function of $\boldsymbol{M}$, with a Lipschitz constant that only depends on distance to $\xi(t)$ of the slit (or slits) nearest to $\xi(t)$. 

Let $\boldsymbol M$ and $\boldsymbol{ M}^*$ be two 		points in moduli space with corresponding standard domains 		$D$ and $D^*$, such that 
\[
	|m_j-m_j^*|,|\theta_j-\theta_j^*|,| \theta_j'-				\theta_j'^*|<\epsilon.
\]
We assume that $\epsilon$ is so small that 
\[
	C_j\cap C_k^*=\emptyset,\text{ whenever }j\neq k.
\]
Denote $z_j, z_j^*$ the endpoints of the slit $C_j$ and $z_j^*, z_j'^*$ the corresponding endpoints of $C_j^*$. Denote $\Psi$ the vector field for $D$ and $\Psi^*$ the vector field for $D^*$. Then we need to show that 
\begin{equation}\label{E:bigO}
	\Psi^*(z_j^*)-\Psi(z_j), \Psi^*(z_j'^*)-\Psi(z_j')
	=O(\epsilon),\quad j=1,\dots,n-1.
\end{equation}

We begin by applying a particular interior variation as in Garabedian's proof of Hadamard's variational formula, \cite{garabedian:1964}.

Denote $x_1, x_2$ real coordinates for $D$. It is easy to see that we can map $D$ one-to-one onto $D^*$ by a transformation
\begin{equation}\label{E:transform}
	x_j^*=x_j+\epsilon S_j,\quad j=1,2,
\end{equation}
which sends the endpoints $z_j, z_j'$ of the slits in $D$ to the endpoints $z_j^*, z_j'^*$, respectively, of the slits in $D^*$,  and 
where 
\[
	S_j=S_j(x_1,x_2),\quad j=1,2,
\]
is a pair of smooth functions in some neighborhood of the closure of $D$. In complex notation, 
\[
	z^*=z+\epsilon F(z,\bar{z}),
\]
where
\[
	F(z,\bar{z})=S_1(x_1,x_2)+i S_2(x_1,x_2).
\] 
Similarly, we put $\zeta^*=\zeta+\epsilon F(\zeta,\bar{\zeta})$, $w^*=w+\epsilon F(w,\bar{w})$. Denote
\[
	g(z,\zeta;\epsilon)=G^*(z^*,w^*)
\]
the transformed Green function. Here $G^*$ is the Green function of $D^*$. Then it is straightforward to show that  $g$ satisfies the equation $L_{\epsilon}[g]=0$, where $L_{\epsilon}$  represents the self-adjoint differential operator defined by
\[
	L_{\epsilon}[g]=\sum_{k,l=1}^2\frac{\partial}{\partial x_k}\left(
	A_{kl} \frac{\partial g}{\partial x_l}\right),
\]
with coefficients
\[
	A_{kl}=\left(\frac{\partial(x_1^*,x_2^*)}{\partial(x_1,x_2)}\right)^{-1}
	\sum_{j=1}^2\frac{\partial x_j^*}{\partial x_k}
	\frac{\partial x_j^*}{\partial x_l}.
\]
Note that $A_{11} A_{22}-A_{12}A_{21}=1$ and that 
\begin{equation}\label{E:A bound}
	A_{kl}=\delta_{kl}+\epsilon\phi_{kl}(x_1,x_2,\epsilon),
\end{equation}
where $\phi_{kl}$ is smooth for  $x_1,x_2$ in a neighborhood of the closure of $D$ and $\epsilon$ in a neighborhood of zero. Let $a_{jk}=a_{jk}(z)$ be the inverse matrix of $A$, and denote $\Gamma$ the quadratic form
\[
	\Gamma(z,\zeta)=\sum_{j,k=1}^2 a_{jk}(\zeta)(x_j-\xi_j)(x_k-\xi_k),
\]
where $\xi_1$ and $\xi_2$ stand for the real and imaginary parts of $\zeta$. Finally, denote $\alpha=\alpha(z,\zeta)$ a fixed smooth function of the two points $z$ and $\zeta$ in $D$, which fulfills the boundary condition $\alpha(z,\zeta)=0$ when either $z$ or $\zeta$ lies on $\partial D$, but has the value
\[
	\alpha(\zeta,\zeta)=\frac{1}{4\pi}
\]
when $z$ coincides with $\zeta$ inside $D$. Then
\[
	P_{\epsilon}(z,\zeta)\equiv\alpha(z,\zeta)
	\ln\frac{1}{\Gamma(z,\zeta)}
\]
defines a parametrix for $L_{\epsilon}$. Direct computation shows that 
\begin{align}\label{E:parametrix}
	L_{\epsilon}[P_{\epsilon}]-L_0[P_0]
	&=\epsilon\phi(z,\zeta,\epsilon),\notag\\
	L_{\epsilon}[G-P_0]-L_0[G-P_0]&=\epsilon\psi(z,\zeta,\epsilon),
\end{align}
where $\phi$ and $\psi$ are smooth functions in $z,\zeta$, and $\epsilon$, except for a $1/r$-singularity on the diagonal (here $r=|z-\zeta|$). Using Green's identity and the boundary conditions, 
\begin{align}\label{E:Green}
	g(\zeta,w;\epsilon)&=G(w,\zeta)-P_0(w,\zeta)
	+P_{\epsilon}(w,\zeta)\notag\\
	&\quad+\int_D g(z,w;\epsilon)L_{\epsilon}[G-P_0+P_{\epsilon}]
	\ dx_1\ dx_2,
\end{align}
see \cite{garabedian:1964}. We note that the integral exists in the sense of Lebesgue as the $1/r$ singularity is integrable. It then follows from \eqref{E:A bound} and \eqref{E:parametrix} that 
\[
	g(\zeta,w;\epsilon)-G(w,\zeta)=O(\epsilon),
\]
uniformly in $\zeta,w$ in the closure of $D$ if $|\zeta-w|$ is bounded away from zero. If we take the normal derivative of both sides of \eqref{E:Green} with respect to $w\in\partial D$, then the identity continues to hold as the $1/|z-w|$-singularity that now appears in the integrand is still integrable. We find 
\[
	\frac{\partial}{\partial n_w}g(\zeta,w;\epsilon)
	-\frac{\partial}{\partial n_w}G(w,\zeta)=O(\epsilon),
\]
uniformly in $\zeta\in\bar{D}$ if $|\zeta-w|$ is bounded away from zero. If we take a further derivative in \eqref{E:Green}, 
\[
	\frac{\partial}{\partial \zeta}=\frac{1}{2}\left(\frac{\partial}{\partial \xi_1}
	-i\frac{\partial}{\partial \xi_2} \right),
\]
then the integrand is no longer integrable in the sense of Lebesgue and has to be understood in the sense of a Cauchy principle value. Thus, \eqref{E:Green}, together with \eqref{E:A bound} and 
\eqref{E:parametrix} still imply
\begin{equation}\label{E:main}
	\frac{\partial^2}{\partial \zeta\partial n_w}g(\zeta,w;\epsilon)
	-\frac{\partial^2}{\partial\zeta\partial n_w}G(w,\zeta)=O(\epsilon),
\end{equation}
uniformly in $\zeta\in\bar{D}$ if $|\zeta-w|$ is bounded away from zero. We now note that all domain functions that are used in the construction of the vectorfield $\Psi$ can be obtained via integration from $\partial^2/(\partial\zeta\partial n_w) G(w,\zeta)$. Thus we will be done if we can show that 
\[
	\frac{\partial^2}{\partial \zeta^*\partial n_{w^*}}G^*(\zeta^*,w^*)
	-\frac{\partial^2}{\partial\zeta\partial n_w}G(w,\zeta)=O(\epsilon),
\]
uniformly in $\zeta\in\bar{D}$ if $|\zeta-w|$ is bounded away from zero. But this follows from \eqref{E:main} and 
\[
	\frac{\partial^2}{\partial \zeta^*\partial n_{w^*}}G^*(\zeta^*,w^*)
	-\frac{\partial^2}{\partial \zeta\partial n_w}g(\zeta,w;\epsilon)
	=O(\epsilon),
\]
the latter being a consequence of \eqref{E:transform}.

The second part of the theorem now follows from general results about ordinary differential equations, exactly as in the simply connected case.
\end{proof}


\section{Radial SLE in multiply connected domains}\label{S:rSLE}

The purpose of this paper is 1) to give a ``natural" construction of conformally invariant measures on ``simple curves" in multiply connected domains, and 2) to study some of the properties of these random curves. We will now motivate, using informal arguments, our particular construction of conformally invariant measures on simple curves. The arguments lead to a small class of processes which contains  radial $SLE_{\kappa}$ in multiply connected domains.    

For a domain $D$ with $n$ non-degenerate boundary continua, a boundary point $z\in\partial D$, and an interior point $w\in D$, let  $W(D,z,w)$ be the set of Jordan arcs in $D$ with endpoints $z$ and $w$. Denote $\{{\mathcal L}_{D,z,w}^{\boldsymbol{M}}\}_{D,z,w}$ a family of probability measures on Jordan arcs in the complex plane such that 
\[
	{\mathcal L}_{D,z,w}^{\boldsymbol{M}}(W(D,z,w))=1,
\]
and where $\boldsymbol{M}=M(D)$.
Such families arise, or are conjectured to arise, as distributions of interfaces of statistical mechanical systems at criticality. Based on these models, e.g. percolation, one expects that the distributions describing the interfaces in different domains with different marked points are related by a Markovian-type property and conformal invariance. Denote $\gamma$ a random Jordan arc with law $\mathcal L_{D,z,w}^{\boldsymbol{M}}$. The Markovian-type property says that if $\gamma'$ is a sub-arc of $\gamma$ which has $z$ as one endpoint and whose other endpoint we denote by $z'$, and if $\boldsymbol{M}'=M(D\backslash\gamma')$, then the conditional law of $\gamma$ given $\gamma'$ is 
\begin{equation}\label{E:markov}
	\text{law}(\gamma|\gamma')=\mathcal L_{D\backslash\gamma',z',w}^{\boldsymbol{M}'}.
\end{equation}
Conformal invariance means that if $f:D\to D'$ is conformal, $z'=f(z)$, $w'=f(w)$, then
\begin{equation}\label{E:ci}
	\mathcal L_{D',z',w'}^{\boldsymbol{M}}
	=f_{*}\mathcal L_{D,z,w}^{\boldsymbol{M}}.
\end{equation}
If \eqref{E:ci} holds, then to understand the family $\{ \mathcal L_{D,z,w}^{\boldsymbol{M}}\}$ it is enough to consider standard domains $D$, take $w=0$, $z=1$, and, by the identification of standard domains with their moduli, we may write 
\[
	\mathcal L_{D,1,0}^{\boldsymbol{M}}
	=\mathcal L^{\boldsymbol{M}}.
\]
In this case there is a natural parametrization of the Jordan arcs we consider. Let 
\[
	s\in[0,\infty)\mapsto\gamma(s)\in\overline{D}
\]
be a Jordan arc in a standard domain $D$ such that 
\[
	\gamma(0)\in S^1,\ \gamma(0,\infty)\subset D\backslash\{0\},
	\text{ and }\lim_{t\to\infty}\gamma(t)=0.
\]
Denote $\boldsymbol{M}=M(D)$ the point in the moduli space corresponding to $D$ and let $g_t^{\boldsymbol M}$ be the canonical mapping from $D\backslash\gamma[0,t]$ onto a standard domain $D_t:=g_t^{\boldsymbol M}(D\backslash\gamma[0,t])$. Then
\[
	\left(g_{t^*}^{\boldsymbol M}\right)'(0)
	<\left(g_t^{\boldsymbol M}\right)'(0),
\]
whenever $t^*<t$, and $\lim_{t\to\infty}\left(g_t^{\boldsymbol M}\right)'(0)=\infty$. Thus we may and always will assume that $\gamma$ is parametrized by the conformal radius of $S^1\cup\gamma(0,\cdot]$ in $0$ with respect to $\partial D$, i.e. so that $\left(g_t^{\boldsymbol M}\right)'(0)=e^t$. This parametrization is natural in the following sense. If $t\ge0$, $\boldsymbol{M}(t)=M(D_t)$, and $\tilde{\gamma}$ is the curve defined by 
\[
	s\in[0,\infty)\mapsto\tilde{\gamma}(s)
	=g_t^{\boldsymbol{M}}(\gamma(t+s)),
\]
then the canonical mapping $g_s^{\boldsymbol{M}(t)}$ from $D_t\backslash\tilde{\gamma}[0,s]$ is given by
\[
	g_s^{\boldsymbol{M}(t)}=g_{t+s}^{\boldsymbol{M}}\circ
	\left(g_t^{\boldsymbol{M}}\right)^{-1},
\]
and so $g_s^{\boldsymbol{M}(t)}(D_t\backslash\tilde{\gamma}[0,s])=D_{t+s}$. In particular 
\[
	\left(g_s^{\boldsymbol{M}(t)}\right)'(0)=e^{t+s} e^{-t}=e^s,
\]
i.e. $\tilde{\gamma}$ is also parametrized by conformal radius. 

Let now $\{g_s^{\boldsymbol{M}}:s\ge0\}$ be the random family of canonical maps corresponding to the random Jordan arcs $\{\gamma[0,s]:s\ge0\}$ in a standard domain $D$, and denote
\[
	\mathcal{L}^{\boldsymbol{M}}
	=\text{law}(\{g_s^{\boldsymbol{M}}:s\ge0\}).
\]
Then, applying first the Markovian-type property and then conformal invariance, \eqref{E:markov}, \eqref{E:ci}, we find
\[
	\text{law}(\{g_{t+s}^{\boldsymbol{M}}:s\ge0\}|
	g_t^{\boldsymbol M})
	=\left(g_t^{\boldsymbol M}\right)_*^{-1}
	\mathcal L^{\boldsymbol{M}(t)}. 
\]
Equivalently,
\begin{equation}\label{E:law1}
	\text{law}\left(\{g_{t+s}^{\boldsymbol{M}}\circ
	\left(g_t^{\boldsymbol M}\right)^{-1}:s\ge0\}|
	g_t^{\boldsymbol M}\right)
	=\text{law}(\{g_{s}^{\boldsymbol{M}(t)}:s\ge0\}).
\end{equation}
By the radial Komatu-Loewner equation, \eqref{E:RKL}, for each $t\ge0$, the $\sigma$-field generated by $g_t^{\boldsymbol M}$ is equal to $\sigma((\theta(r),\boldsymbol{M}(r)):r\in[0,t])$, where $\exp(i\theta(r))=\xi(r)$, and $\theta(0)=0$. Similarly, it is easy to see that we can reconstruct $g_{t+s}^{\boldsymbol{M}}\circ
	\left(g_t^{\boldsymbol M}\right)^{-1}$ from $\{(\theta(t+r)-\theta(t),\boldsymbol{M}(t+r)):r\in[0,s]\}$. Thus \eqref{E:law1} implies
\begin{align}\label{E:m}
	\text{law}&(\{(\theta(t+s)-\theta(t),\boldsymbol{M}(t+s)):s\ge0\}|
	\{(\theta(r),\boldsymbol{M}(r)):r\in[0,t]\})\notag\\
	&=\text{law}(\{(\tilde{\theta}(s),\tilde{\boldsymbol{M}}(s)):s\ge0\}),
\end{align}
where $\tilde{\boldsymbol{M}}(s)=M(D_t\backslash\tilde{\gamma}[0,s])$, for a random Jordan arc $\tilde{\gamma}$ with law $\mathcal L^{\boldsymbol{M}(t)}$. The equality \eqref{E:m} is precisely the statement that $\{(\theta(t),\boldsymbol{M}(t)):t\ge0\}$ is a Markov process. We note that in the simply connected case ($n=1$), \eqref{E:m} reduces to 
\[
	 \text{law}(\{\theta(t+s)-\theta(t):s\ge0\}|
	\{\theta(r):r\in[0,t]\})
	=\text{law}(\{\tilde{\theta}(s):s\ge0\}),
\]
from which it follows that $\theta$ is a process with independent, and identically distributed increments. From this, continuity, and the symmetry $\text{law}(\theta)=\text{law}(-\theta)$, Schramm derived in \cite{schramm:2000} that $\theta(t)=\sqrt{\kappa}B_t$ for a standard one-dimensional Brownian motion and a positive constant $\kappa$. The continuity follows from the continuity of the Jordan arcs, and the symmetry is actually observed in various discrete models, such as loop-erased random walk.

In the multiply connected case, we also have continuity. Let us now study the It\^o differential for the diffusion $(\theta(t),\boldsymbol{M}(t))$. By Lemma \ref{L:mm}, $\boldsymbol{M}(t)$ is a finite variation process and its differential is given by \eqref{E:mm}. On the other hand, we know from \cite{rohde.schramm:2003}
 that the qualitative properties of $\gamma(t)$ would change with $t$ if the martingale part of $\theta(t)$ has quadratic variation which is nonlinear in $t$. Thus $d\theta(t)=\sqrt{\kappa}dB_t+\text{``drift"}$, and the only open question concerning the diffusion $(\theta(t),\boldsymbol{M}(t))$ is the drift of $\theta(t)$. For a general drift, which may be a function of $\boldsymbol M$, $\theta$, and $\kappa$, the resulting family $\{g_t^D\}$ is a {\it random Loewner chain} and we call the diffusion a {\it Schiffer diffusion}. In the case of percolation, the drift is easily identified. Consider a honeycomb lattice-approximation to our domain $D$. The beginning of the random simple curve $\gamma$ is an edge $e$ of the lattice in $D$, which is horizontal and whose one endpoint is the point 1. The next step is  either up or down (with slope $2\pi/3$). For percolation, each of these possibilities has probability 1/2. Denote $p$ the endpoint in $D$ of the  slit made up of two edges, if the second step went up, and $q$ the corresponding endpoint, if the second step went down. Denote $g$ the canonical map from $D\backslash\{e\}$.  
To derive the drift from this condition, we compare the images of the endpoints $p$ and $q$. Recall that for the unit disk $\mathbb D$
\[
	F(z,w)=\ln(1-\overline{z}w)-\ln(w-z).
\]
Hence
\[
	-\frac{\partial F(z,w)}{\partial n_z}=\frac{z+w}{z-w}, \quad z\in S^1.
\]
For a standard domain $D$, let
\[
	k_D(z,w)=-\frac{\partial F_D(z,w)}{\partial n_z}
	-\frac{z+w}{z-w},\quad z\in S^1,
\]
and define
\[
	k_D(z)=\lim_{w\to z}k_D(z,w),\quad z\in S^1.
\]

If the lattice size is small, say a single edge is length $\sqrt{\epsilon}$, then $g'(0)=\epsilon+o(\epsilon)$. Thus
\begin{align}
	\frac{\Re(g(p)+g(q))}{2}&=\epsilon\lim_{x\to0}\frac{\Psi_D(1-x,1)+\Psi_D(1+x,1)}{2}+o(\epsilon)\notag\\
	&=\epsilon(k_D(1)+\boldsymbol R_D(1)\boldsymbol{P}_D^{-1}\frac{\partial\boldsymbol\omega_D(1)}{\partial n})+o(\epsilon).
\end{align}

Hence, to model cluster-boundaries of percolation in a multiply connected domain $D$ we make the ansatz
\begin{align}\label{E:ansatz}
	d\theta(t)=&-i\left( k(\xi(t);\boldsymbol{M}(t))+{\bf R}(\xi(t);\boldsymbol{M}(t))^T{\bf P}_{\boldsymbol{M}(t)}^{-1}\frac{\partial\boldsymbol\omega(\xi(t);\boldsymbol{M}(t))}{\partial n}\right) dt\notag\\
	&+\sqrt{\kappa}\ dB_t,
\end{align}
with $\xi(t)=e^{i\theta(t)}$, and where $\boldsymbol{M}(t)$ satisfies \eqref{E:mm}.

For other discrete models the same reasoning as above would lead to different drifts. For example, for loop-erased random walk, the probability of stepping up is not the same as stepping down, depending on the configuration of the concentric circular slits. However, the respective probabilities can be calculated in terms of harmonic measure. This leads to a different Schiffer diffusion and we expect the resulting family $\{g_t^D\}$ to be closely related to the  ``harmonic random Loewner chains'' studied by Zhan in his thesis, \cite{zhan:2004}. We leave the question of which Schiffer diffusion corresponds to which discrete model to a forthcoming paper. In principle, representation-theoretic considerations as are done in conformal field theory should identify the relevant class of Schiffer diffusions.
 

\section{Locality}\label{S:locality}

In this section we show that the ansatz \eqref{E:ansatz} leads to random growing compacts satisfying the locality property if $\kappa=6$. 
Denote $\{g_t^E,t\ge0\}$ the solution of the radial Komatu-Loewner equation in a standard domain $E$  starting at $z=1$ for the Schiffer diffusion \eqref{E:ansatz}. Denote  $\{K_t,t\ge0\}$ the associated growing compacts. Let $A$ be a hull in $E$ that contains neither zero nor $z=1$. For the following calculations we restrict to the event $\{t<\tau\}$, where $\tau:=\inf\{t:K_t\cap A\neq\emptyset\}$. Let $\Phi_A$ be the canonical mapping from $E\backslash A$, $g_t^*$ the canonical mapping from $\Phi_A(E\backslash(A\cup K_t))$, and $h_t$ the canonical mapping from  $g_t(E\backslash(A\cup K_t))$. Since the canonical mapping for $E\backslash(A\cup K_t)$ is unique, we have
\begin{equation}
	h_t\circ g_t=g_t^*\circ\Phi_A.
\end{equation}
Furthermore, up to a time change, the family $\{g_t^*\}$ also satisfies a radial Komatu-Loewner equation beginning with the standard domain $E^*:=\Phi_A(E\backslash A)$. In fact, it follows from \eqref{E:change of r} that 
\begin{align}\label{E:KLstar}
	\partial_t&\ln g_t^*(z)\notag\\
	&=|h_t'(\xi(t))|^2\left(\frac{\partial F^*(\xi^*(t),w_t^*;t)}{\partial n}+{\bf R}^*(w_t^*;t)^T\left({\bf P}_t^*\right)^{-1}\frac{\partial\boldsymbol\omega^*(\xi^*(t);t)}{\partial n}\right),
\end{align}
where $w_t^*=g_t^*(z)$, and $\xi^*(t)=h_t(\xi(t))$. The question we are interested in is whether $(\xi^*,\boldsymbol{M}^*)$ is a time change of $(\xi,\boldsymbol{M})$. Since $h_t=g_t^*\circ\Phi_A\circ g_t^{-1}$, we have
\begin{equation}\label{E:ht}
	\partial_t h_t(z)=\left[\partial_t g_t^*\right](\Phi_A(g_t^{-1}(z)))+(g_t^*\circ\Phi_A)'(g_t^{-1}(z))(\partial_t g_t^{-1}(z)),
\end{equation}
and we note that 
\begin{equation}\label{E:g inverse}
	\partial_t g_t^{-1}(z)=-(g_t^{-1})'(z)z\left(\frac{\partial F(\xi(t),z;t)}{\partial n}+{\bf R}(z;t)^T{\bf P}_t^{-1}\frac{\partial\boldsymbol\omega(\xi(t);t)}{\partial n}\right).
\end{equation}
Let $\phi_t(z)=-i\ln h_t(e^{iz})$. Then $\partial_t \phi_t(z)=-i\partial_t h_t(e^{iz})/h_t(e^{iz})$, and $\phi_t'(z)=e^{iz}h_t'(e^{iz})/h_t(e^{iz})$.
Thus \eqref{E:ht},\eqref{E:KLstar}, and \eqref{E:g inverse} imply
\begin{align}\label{E:phit}
	\partial_t\phi_t(z)&=\frac{\phi_t'(\theta(t))^2}{i}\frac{\partial F^*(\xi^*(t),e^{i\phi_t(z)};t)}{\partial n}\notag\\
	&\quad+\frac{\phi_t'(\theta(t))^2}{i}{\bf R}^*(e^{i\phi_t(z)};t)^T\left({\bf P}_t^*\right)^{-1}\frac{\partial\boldsymbol\omega^*(\xi^*(t);t)}{\partial n}\notag\\
	&\quad-\frac{\phi_t'(z)}{i}\left(\frac{\partial F(\xi(t),e^{iz};t)}{\partial n}+{\bf R}(e^{iz};t)^T{\bf P}_t^{-1}\frac{\partial\boldsymbol\omega(\xi(t);t)}{\partial n}\right).
\end{align}
Hence the stochastic differential 
\[
	\partial_t\phi_t(z)\ dt+\phi_t'(\theta(t))\ d\theta(t)
\]
has martingale part $\phi_t'(\theta(t))\sqrt{\kappa}\ dB_t$ and its drift part consists of the three components
\begin{align}
	I:&=\frac{\phi_t'(\theta(t))^2}{i}\left(\frac{\partial F^*(\xi^*(t),e^{i\phi_t(z)};t)}{\partial n}-k^*(\xi^*(t);t)\right) dt\notag\\
	&\quad-\frac{\phi_t'(z)}{i}\left(\frac{\partial F(\xi(t),e^{iz};t)}{\partial n}-k(\xi(t);t)\right) dt,\notag\\
	II:&=\frac{\phi_t'(\theta(t))^2}{i}\left(k^*(\xi^*(t);t)+{\bf R}^*(e^{i\phi_t(z)};t)^T({\bf P}_t^*)^{-1}\frac{\partial\boldsymbol\omega^*(\xi^*(t);t)}{\partial n}\right)dt,\notag\\
	III:&=\frac{\phi_t'(\theta(t))-\phi_t'(z)}{i}\left(k(\xi(t);t)+{\bf R}(e^{iz};t)^T{\bf P}_t^{-1}\frac{\partial\boldsymbol\omega(\xi(t);t)}{\partial n}\right) dt.
\end{align}
When $z\to\theta(t)$, then part $III$ converges to zero, and part $II$, together with the martingale part, converges to a time-change  of \eqref{E:ansatz} starting at $E^*$. Finally, the limit of part $I$ is by the definition of $k(\xi;t)$ equal to 
\begin{equation}
	\lim_{z\to \theta}\left(\frac{2\phi'(\theta)^2}{\phi(z)-\phi(\theta)}-
	\frac{2\phi'(z)}{z-\theta}\right)=-3\phi''(\theta).
\end{equation}
Thus, by It\^o's formula,
\begin{align}\label{E:ito}
	d\phi_t&(\theta(t))\notag\\
	&=\frac{\phi_t'(\theta(t))^2}{i}\left(k^*(\xi^*(t);t)+{\bf R}^*(\xi^*(t);t)^T({\bf P}_t^*)^{-1}\frac{\partial\boldsymbol\omega^*(\xi^*(t);t)}{\partial n}\right) dt\notag\\
	&\quad+\frac{\kappa-6}{2}\phi_t''(\theta(t))\ dt+\phi_t'(\theta(t))\sqrt{\kappa}\ dB_t,
\end{align}
which is indeed a time-change of \eqref{E:ansatz} if and only if $\kappa=6$. From \eqref{E:KLstar} it follows immediately that the equations for $\boldsymbol{M}^*$ are given by the same time change of the equations for $\boldsymbol{M}$.  

\begin{theorem}[Radial $\text{SLE}_6$]
	The solution to the radial Komatu-Loewner equation based on the Schiffer diffusion \eqref{E:ansatz} satisfies the locality property if and only if $\kappa=6$.
\end{theorem}

\end{document}